\documentstyle[twoside]{article}
\title{\bf Uniform asymptotics of a Gauss hypergeometric function with two large parameters, V}
\author{\sc R. B. Paris\footnote{E-mail address:\ \ {\tt r.paris@abertay.ac.uk}}\\
\\
{\em Division of Computing and Mathematics,}\\
{\em Abertay University, Dundee DD1 1HG, UK}\\
}

\begin{document}
\newcommand{\bee}{\begin{equation}}
\newcommand{\ee}{\end{equation}}
\def\f#1#2{\mbox{${\textstyle \frac{#1}{#2}}$}}
\def\dfrac#1#2{\displaystyle{\frac{#1}{#2}}}
\newcommand{\fr}{\frac{1}{2}}
\newcommand{\fs}{\f{1}{2}}
\newcommand{\g}{\Gamma}
\newcommand{\la}{\lambda}
\newcommand{\al}{\alpha}
\newcommand{\br}{\biggr}
\newcommand{\bl}{\biggl}
\renewcommand{\topfraction}{0.9}
\renewcommand{\bottomfraction}{0.9}
\renewcommand{\textfraction}{0.05}
\newcommand{\gtwid}{\raisebox{-.8ex}{\mbox{$\stackrel{\textstyle >}{\sim}$}}}
\newcommand{\ltwid}{\raisebox{-.8ex}{\mbox{$\stackrel{\textstyle <}{\sim}$}}}
\newcommand{\mcol}{\multicolumn}
\date{}
\maketitle
\pagestyle{myheadings}
\markboth{\hfill {\it R.B. Paris} \hfill}
{\hfill {\it Hypergeometric function with large parameters} \hfill}
\begin{abstract} 
We consider the uniform asymptotic expansion for the Gauss hypergeometric function
\[{}_2F_1(a+\epsilon\lambda,b;c+\lambda;x),\qquad 0<x<1\]
as $\lambda\to+\infty$ in the neigbourhood of $\epsilon x=1$ when the parameter $\epsilon>1$ and the constants $a$, $b$ and $c$ are supposed finite. Use of a standard integral representation shows that the problem reduces to consideration of a simple saddle point near an endpoint of the integration path. A uniform asymptotic expansion is first obtained by employing Bleistein's method. An alternative form of uniform expansion is derived following the approach described in Olver's book [{\it Asymptotics and Special Functions}, p.~346]. This second form has several advantages over the Bleistein form.

Numerical results illustrating the accuracy of the different expansions are given. 
\vspace{0.4cm}

\noindent {\bf MSC:} 33C05, 34E05, 41A60
\vspace{0.3cm}

\noindent {\bf Keywords:} Hypergeometric function, uniform asymptotic expansion, large parameters, saddle near an endpoint\\
\end{abstract}

\vspace{0.2cm}

\noindent $\,$\hrulefill $\,$

\vspace{0.2cm}

\begin{center}
{\bf 1. \  Introduction}
\end{center}
\setcounter{section}{1}
\setcounter{equation}{0}
\renewcommand{\theequation}{\arabic{section}.\arabic{equation}}
Asymptotic expansions for large $\lambda$ of the Gauss hypergeometric function
\[{}_2F_1(a+\epsilon_1\la,b+\epsilon_2\la;c+\epsilon_3\la;z)\]
for $\epsilon_r=\pm 1$ ($1\leq r\leq 3$), finite values of the parameters $a$, $b$ and $c$ and fixed complex $z$ were first considered by Watson \cite{W} in 1918. More recently, this study was extended in \cite{P1, P2} for $\epsilon_r>0$, and also in \cite{Cvit} for the case of two large parameters.  

In \cite[\S 3]{P1}, the case when $\epsilon_2=0$ was considered which, after a rescaling of the parameter $\la$ and setting $\epsilon=\epsilon_1/\epsilon_3$, can be expressed as the function 
\bee\label{e11}
{}_2F_1(a+\epsilon\la,b;c+\la;z),\qquad \epsilon>0.
\ee
The method of steepest descents was applied to suitable integral representations valid for $0<\epsilon<1$ and $\epsilon>1$.
The resulting expansion when $\epsilon<1$ is given in \cite[(3.6)]{P1} and has the leading behaviour $(1-\epsilon z)^{-b}$; this case is also considered in \cite{LP}.
The asymptotic expansion of (\ref{e11}) when $\epsilon>1$ is given by \cite[(3.13)]{P1}
\bee\label{e12}
\frac{{\hat G}(\la)}{\sqrt{2\pi} (1-\epsilon z)^b} \,\epsilon^{a-\frac{1}{2}+\epsilon\lambda}(\epsilon-1)^{c-a-\frac{1}{2}+(1-\epsilon)\lambda}\, \sum_{k=0}^\infty \frac{c_k \g(k+\fs)}{\lambda^{k+\frac{1}{2}} \g(\fs)}
\ee
as $\lambda\to\infty$, where 
\[{\hat G}(\la)=\frac{\g(c+\lambda)\g(1+a-c+(\epsilon-1)\lambda)}{\g(a+\epsilon\lambda)},\]
$c_0=1$ and the coefficients $c_1$ and $c_2$ are explicitly stated in \cite[(3.4), (3.5)]{P1}. Application of Stirling's formula
$\g(a+x)\sim \sqrt{2\pi}\,x^{x+a-\frac{1}{2}} e^{-x}$ $(x\to+\infty)$
shows that
\[{\hat G}(\la)\sim (2\pi\lambda)^\frac{1}{2} \epsilon^{\frac{1}{2}-a-\epsilon\lambda} (\epsilon-1)^{\frac{1}{2}+a-c+(\epsilon-1)\lambda}\qquad(\lambda\to+\infty),
\]
so that the leading large-$\lambda$ behaviour of (\ref{e11}) again reduces to $(1-\epsilon z)^{-b}$; see also \cite{Cvit}.

The above Poincar\'e-type expansion is subject to an inconvenient restriction (when $\epsilon>1$), which results from the requirement that in the integral used to determine its expansion the singularity of the integrand at $1/z$ should lie outside the closed-loop contour. In the case of real $z$ ($=x$) with $0<x<1$, which is the situation we consider throughout in the present paper, this restriction corresponds to the condition $\epsilon x<1$.
Thus, the result presented in (\ref{e12}) does not cover the situation when $\epsilon x>1$ nor that in the neighbourhood of $\epsilon x=1$. 

In a recent paper \cite{H}, Harper encountered a case of (1.1) with $\epsilon=\f{3}{2}$ and $x=\f{2}{3}$ (so that $\epsilon x=1$) arising in the theory of a gas bubble rising in a solution of a substance that raises the surface tension of the liquid. This case corresponds to a saddle point coincident with an endpoint of integration, for which an expansion can be found by application of Laplace's method; see \cite[\S 3]{H}. A uniform expansion of (\ref{e11}) valid in the neighbourhood of $\epsilon x=1$ when $b=1$ was considered in \cite{P18} by the standard procedure of the method of steepest descents modified to deal with the situation when a saddle point is near a simple pole. This study was incomplete since the case of integer $b\geq 2$ was dealt with by the use of a recursion relation and did not cover the case of non-integer values of $b$.

Our aim in this paper is to derive the expansion of the hypergeometric function\footnote{Note that when $\epsilon<1$ the function $F(a+\epsilon\la,b;c+\la;z)$ converges at $z=1$, since the convergence condition $ c-a-b+(1-\epsilon)\lambda$ will be positive as $\lambda\to+\infty$. When $\epsilon>1$, however, this condition will be broken for sufficiently large $\lambda$ and the function will not converge at $z=1$; see \cite[\S 1]{P18}.} in (\ref{e11}) when $z=x$ with $0<x<1$, $\epsilon>1$ and real $b>0$ that holds uniformly in the neighbourhood of $\epsilon x=1$.  
The integral representation we choose to represent this function is characterised by a simple saddle point situated near an endpoint of the integration path. We shall apply the standard Bleistein
expansion to deal with the case. This results in an expansion involving the parabolic cylinder function and its derivative as approximating functions, which describe the transition of the saddle through the endpoint. A drawback with this method is the calculation of the higher coefficients in the resulting uniform expansion. An alternative expansion is obtained following the approach described by Olver in \cite[p.~346]{O}, which yields a similar expansion to the Bleistein expansion but which does not present the same inconvenience in the determination of the higher coefficients. 
In a final section, numerical results are presented to demonstrate the accuracy of the two types of expansion obtained.

\vspace{0.6cm}

\begin{center}
{\bf 2. \ The Bleistein expansion of ${}_2F_1(a+\epsilon\la,b;c+\la;x)$ for $\lambda\to+\infty$}
\end{center}
\setcounter{section}{2}
\setcounter{equation}{0}
\renewcommand{\theequation}{\arabic{section}.\arabic{equation}}
We consider the function defined in (\ref{e11}) as $\la\to+\infty$ when $0<x<1$, $\epsilon>1$ and the parameters $a$, $b>0$, $c$ are real and finite. From \cite[(15.6.1)]{DLMF} we have the integral representation
\bee\label{e21}
{}_2F_1(a+\epsilon\la,b;c+\la;x)=\frac{G(\la)}{\g(b)} \int_0^1 t^{b-1} (1-t)^{c-b-1+\la} (1-xt)^{-a-\epsilon\la}dt
\ee
provided $c-b+\la>0$, where \cite[(5.11.3)]{DLMF}
\bee\label{e22}
G(\la):=\frac{\g(c+\la)}{\g(c+\la-b)}\sim \la^b\sum_{k=0}^\infty \frac{h_k(b,c)}{\la^k}
\ee
as $\la\to+\infty$, with the first few coefficients $h_k(b,c)$ given by
\[h_0(b,c)=1,\quad h_1(b,c)=\fs b(2c-b-1),\]
\[ h_2(b,c)=\f{1}{24}b(b-1)\{3(2c-b-1)^2-b-1\}, \ldots\ .\]

The above integral can be cast in the form
\bee\label{e23}
{}_2F_1(a+\epsilon\la,b;c+\la;x)=\frac{G(\la)}{\g(b)}\int_0^1t^{b-1} f(x,t) e^{-\la\psi(x,t)}dt,
\ee
where
\[\psi(x,t)=\epsilon \log\,(1-xt)-\log\,(1-t),\qquad f(x,t)=\frac{(1-t)^{c-b-1}}{(1-xt)^a}.\]
The phase function $\psi(x,t)$ has a saddle point at $\partial \psi(x,t)/\partial t=0$; that is, at the point $t=t_s$ where
\bee\label{e24}
t_s=\frac{\epsilon x-1}{(\epsilon-1)x}.
\ee
When $\epsilon x>1$, the saddle point lies in the interval $(0,1)$ and, since $x$ is assumed to satisfy $0<x<1$, it is easily verified that the steepest descent path through the saddle point is the integration path $0\leq t\leq 1$. When $\epsilon x=1$, $t_s=0$ so that the saddle coincides with the endpoint $t=0$. When $0<\epsilon x<1$ the saddle $t_s<0$, with $t_s\to-\infty$ as $x\to 0+$. 

We now introduce the new variable $u$ by
\bee\label{e24a}
\psi(x,t)-\psi(x,t_s)=\fs(u-\al)^2,
\ee
where the point $u=\al$ corresponds to the saddle $t=t_s$. The points $u=0$ and $t=0$ are made to correspond by the choice
\bee\label{e24b}
\al=\pm (2\{\psi(x,0)-\psi(x,t_s)\})^{1/2}=\pm(-2\psi(x,t_s))^{1/2}
\ee
since $\psi(x,0)=0$,
where the upper or lower sign is chosen according as $t_s>0$ or $t_s<0$, respectively; thus $\al$ {\it has the same sign as} $t_s$. Substitution of the new variable $u$ into (\ref{e23}) then produces
\bee\label{e25}
{}_2F_1(a+\epsilon\la,b;c+\la;x)=\frac{G(\la)}{\g(b)}\,e^{-\la\psi(x,t_s)} \int_0^\infty u^{b-1}g_0(u) e^{-\fr\la(u-\al)^2}du,
\ee
where 
\bee\label{e25a}
g_0(u)\equiv g_0(x,u)=f(x,t) \bl(\frac{t}{u}\br)^{\!\!b-1} \frac{dt}{du}~.
\ee

The function $g_0(u)$ is decomposed into a sequence of functions $g_k(u)$ following the procedure introduced by Bleistein \cite{B}; for a description of the method see, for example, \cite{J}, \cite[p.~360]{Wong} or \cite[p.~62]{PBook}. We write
\bee\label{e25b}
g_k(u)=A_k(\al)+B_k(\al)(u-\al)+u(u-\al) G_k(u)\qquad (k\geq 0)
\ee
with
\bee\label{e25c}
g_{k+1}(u)=bG_k(u)+uG_k'(u).
\ee
The coefficients $A_k(\al)$, $B_k(\al)$ are given by
\[A_k(\al)=g_k(\al),\qquad B_k(\al)=\frac{1}{\al}\{g_k(\al)-g_k(0)\}.\]
Substitution of $g_0(u)$ given by (\ref{e25a}) into (\ref{e25}), combined with an integration by parts, then yields after repetition of this process for $k\geq 1$ the following standard expansion:
\newtheorem{theorem}{Theorem}
\begin{theorem}$\!\!\!.$ Let $\epsilon>1$ and the parameters $a$, $b$ and $c$ be real with $b>0$. Further let $0<x<1$. Then we have the uniform expansion valid in the neighbourhood of $\epsilon x=1$ $(\al=0)$
\bee\label{e26}
F(a+\epsilon\la,b;c+\la;x)\sim G(\la)e^{-\la\psi(x,t_s)}\bl\{\frac{W_b(\al\la^\fr)}{\la^{\fr b}}\sum_{k=0}^\infty \frac{A_k(\al)}{\la^k}+\frac{W_b'(\al\la^\fr)}{\la^{\fr b+\fr}}\sum_{k=0}^\infty \frac{B_k(\al)}{\la^k}\br\}
\ee
as $\la\to+\infty$, where
the quantity $W_b(\chi)$ can be expressed in terms of the parabolic cylinder function $D_\nu(z)$ by
\bee\label{e27}
W_b(\chi)=\frac{1}{\g(b)}\int_0^\infty \tau^{b-1} e^{-\fr(\tau-\chi)^2}d\tau=e^{-\chi^2/4}D_{-b}(-\chi).
\ee
\end{theorem}
\vspace{0.2cm}

The leading coefficients in the above expansion are given by
\bee\label{e28}
A_0(\al)=g_0(\al),\qquad B_0(\al)=\frac{1}{\al}\{g_0(\al)-g_0(0)\},
\ee
where, noting that $f(x,0)=1$,
\bee\label{e28a}
g_0(0)=\bl(\frac{-\al}{\psi'(x,0)}\br)^{\!b},\qquad g_0(\al)=f(x,t_s) \bl(\frac{t_s}{\al}\br)^{\!\!b-1} \{\psi''(x,t_s)\}^{-1/2}
\ee
with primes on $\psi$ denoting differentiation with respect to $t$. The next two coefficients are
\bee\label{e29}
A_1(\al)=g_1(\al),\qquad B_1(\al)=\frac{1}{\al}\{g_1(\al)-g_1(0)\},
\ee
where (see Appendix A for details)
\bee\label{e210}
g_1(0)=\frac{b}{\al}\{B_0(\al)-g_0'(0)\},\qquad g_1(\al)=\frac{(1-b)}{\al}\{B_0(\al)-g'_0(\al)\}+\fs g_0''(\al),
\ee
with $g_0'(0)$ given in (\ref{a6}).
%\bee\label{e211}
%g_0'(0)=g_0(0)\bl\{\frac{(b+1)}{2\al}\bl(\frac{\al^2 \psi''(x,0)}{\psi'(x,0)^2}-1\br)-\frac{\al(ax+b+1-c)}{\psi'(x,0)}\br\}.
%\ee
The computation of $g_0'(\al)$ and $g_0''(\al)$ in specific cases is discussed in Section 6. The higher coefficients $A_k(\al)$ and $B_k(\al)$ are more difficult to compute. 
\vspace{0.6cm}

\begin{center}
{\bf 3. \ An alternative uniform expansion of ${}_2F_1(a+\epsilon\la,b;c+\la;x)$ for $\lambda\to+\infty$}
\end{center}
\setcounter{section}{3}
\setcounter{equation}{0}
\renewcommand{\theequation}{\arabic{section}.\arabic{equation}}
To derive an alternative form of uniform expansion for ${}_2F_1(a+\epsilon\la,b;c+\la;x)$ we follow the procedure described in \cite[p.~346]{O}. The function $g_0(u)$ defined in (\ref{e25a}) is expanded as a Taylor series about the saddle $u=\al$ to yield
\bee\label{e31}
g_0(u)=\sum_{r=0}^\infty p_k(\al) (u-\al)^k,\qquad p_k(\al)=\frac{g_0^{(k)}(\al)}{k!}\qquad (|u|<\infty),
\ee
where the coefficients $p_k(\al)$ are continuous at $\al=0$. Substitution into (\ref{e25}) then yields (formally)
\[{}_2F_1(a+\epsilon\la,b;c+\la;x)\sim \frac{G(\la) e^{-\la\psi(x,t_s)}}{\g(b)} \sum_{k=0}^\infty p_k(\al) \int_0^\infty u^{b-1} (u-\al)^k e^{-\fr\la(u-\al)^2}du\]
\bee\label{e32}
=G(\la) e^{-\la\psi(x,t_s)} \sum_{k=0}^\infty \frac{p_k(\al) S_k(\al\la^\fr)}{\la^{\fr b+\fr k}},
\ee
where
\[S_k(\chi)=\frac{1}{\g(b)} \int_0^\infty \tau^{b-1} (\tau-\chi)^k e^{-\fr (\tau-\chi)^2}d\tau.\]

Upon replacement of the integration variable by $\tau=s+\chi$ in $S_k(\chi)$, it can be shown by Laplace's method that as $\chi\to\+\infty$ (cf. \cite[(9.09), p.~346]{O})
\[S_k(\chi) \sim \left\{\begin{array}{ll} \dfrac{2^{(k+1)/2}}{\g(b)}\,\g(\fs k+\fs) \chi^{b-1} & (k\ \mbox{even})\\
\dfrac{2^{k/2+1}}{\g(b-1)}\,\g(\fs k+1) \chi^{b-2} & (k\ \mbox{odd}). \end{array}\right.\]
Hence for fixed non-zero $\al$ and large $\chi$ successive terms of (\ref{e32}), taken in pairs, are of decreasing asymptotic order. 

It is readily seen that
\bee\label{e32a}
S_0(\chi)=W_b(\chi)=e^{-\chi^2/4}D_{-b}(-\chi),\qquad S_1(\chi)=W_b'(\chi)=e^{-\chi^2/4} D_{-b+1}(-\chi).
\ee
Higher $k$-values of $S_k(\chi)$ satisfy the recurrence relation (cf. \cite[Ex.~(9.1), p.~346]{O})
\bee\label{e33a}
S_k(\chi)=-\chi S_{k-1}(\chi)+(b+k-2) S_{k-2}(\chi)+\chi (k-2) S_{k-3}(\chi)\qquad (k\geq 2).
\ee
This enables us to express $S_k(\chi)$ for $k\geq 2$ in terms of $S_0(\chi)$ and $S_1(\chi)$, and hence in terms of $W_b(\chi)$ and its derivative by (\ref{e32a}), in the form
\bee\label{e33}
S_k(\chi)=c_k(\chi) S_0(\chi)+d_k(\chi) S_1(\chi)\qquad(k\geq 0).
\ee
The coefficients $c_k(\chi)$, $d_k(\chi)$ are listed in Table 1 for $0\leq k\leq 7$.
\begin{table}[th]
\caption{\footnotesize{Values of the coefficients $c_k(\chi)$ and $d_k(\chi)$ for $k\leq 7$ where $\mu:=b-1$}} \label{t1}
\begin{center}
\begin{tabular}{|l|l|l|}
\hline
&&\\[-0.3cm]
\mcol{1}{|c|}{$k$} & \mcol{1}{c|}{$c_k(\chi)$} & \mcol{1}{c|}{$d_k(\chi)$}\\
\hline
&&\\[-0.3cm]
0 & 1 & 0\\
1 & 0 & 1\\
2 & $\mu\!+\!1$ & $-\chi$\\
3 & $-\mu\chi$& $\mu\!+\!2\!+\!\chi^2$\\
4 & $(\mu\!+\!1)(\mu\!+3)\!+\!\mu\chi^2$ & $-(3\!+\!2\mu)\chi\!-\!\chi^3$\\
5 & $-(5\!+\!2\mu)\mu\chi\!-\!\mu\chi^3$ & $(\mu\!+\!2)(\mu\!+\!4)\!+\!(4\!+\!3\mu)\chi^2\!+\!\chi^4$\\
6 & $(\mu\!+\!1)(\mu\!+\!3)(\mu\!+\!5)\!+\!(6\!+\!3\mu)\mu\chi^2\!+\!\mu\chi^4$ & $-(15\!+\!15\mu\!+\!3\mu^2)\chi\!-\!(5\!+\!4\mu)\chi^3\!-\!\chi^5$\\
7 & $-(33\!+\!21\mu\!+\!3\mu^2)\mu\chi\!-\!(7\!+\!4\mu)\mu\chi^3\!-\!\mu\chi^5$ &   $(\mu\!+\!2)(\mu\!+\!4)(\mu\!+\!6)\!+\!(24\!+\!27\mu\!+\!6\mu^2)\chi^2$\\
&& $\hspace{3cm}+(6\!+\!5\mu)\chi^4\!+\!\chi^6$ \\
[.1cm]\hline\end{tabular}
\end{center}
\end{table}

Then if we define the coefficients
\bee\label{e34}
C_k(\al,\la)=p_k(\al) c_k(a\la^\fr),\qquad D_k(\al,\la)=p_{k+1}(\al) d_{k+1}(a\la^\fr)\qquad (k\geq0),
\ee
the expansion appearing in (\ref{e32}) becomes
\[\sum_{k=0}^\infty \frac{p_k(\al) S_k(\al\la^\fr)}{\la^{b/2+k/2}}=\frac{1}{\la^{\fr b}}\sum_{k=0}^\infty\bl\{\frac{C_k(\al,\la)S_0(\al\la^\fr)}{\la^{\fr k}}+\frac{D_k(\al,\la)S_1(\al\la^\fr)}{\la^{\fr k+\fr}}\br\}.\]
After some routine algebra using the values of the coefficients listed in Table 1, we find the expansion {\it in inverse integer powers of $\la$} given by
\[\sum_{k=0}^\infty \frac{p_k(\al) S_k(\al\la^\fr)}{\la^{\fr b+\fr k}}=\frac{S_0(\al\la^\fr)}{\la^{\fr b}}\sum_{k=0}^\infty\frac{{\cal C}_k(\al)}{\la^k}+\frac{S_1(\al\la^\fr)}{\la^{\fr b+\fr}}\sum_{k=0}^\infty\frac{{\cal D}_k(\al)}{\la^k},\]
where the first few coefficients are 
\begin{eqnarray}
{\cal C}_0(\al)&=&p_0(\al),\nonumber\\
{\cal C}_1(\al)&=&p_2(\al)+\mu(p_2(\al)-\al p_3(\al)+\al^2 p_4(\al)-\al^4 p_5(\al)+\al^4 p_6(\al)+\cdots),\nonumber\\
{\cal C}_2(\al)&=&(\mu+1)(\mu+3)p_4(\al)-\mu((5+2\mu)\al p_5(\al)+(6+3\mu)\al^2 p_6(\al)+\cdots)\,,\nonumber\\
{\cal C}_3(\al)&=&(\mu+1)(\mu+3)(\mu+5)p_6(\al)+\cdots \label{e34a}
\end{eqnarray}
and
\begin{eqnarray}
{\cal D}_0(\al)&=&p_1(\al)-\al p_2(\al)+\al^2 p_3(\al)-\al^3 p_4(\al)+\al^4 p_5(\al)-\al^5 p_6(\al)+\al^6 p_7(\al)- \cdots\, ,\nonumber\\
{\cal D}_1(\al)&=&(2+\mu)p_3(\al)-(3+2\mu)\al p_4(\al)+(4+3\mu)\al^2p_5(\al)-(5+4\mu)\al^3 p_6(\al)\nonumber\\
&& +(6+5\mu)\al^4 p_7(\al)+\cdots,\nonumber\\
{\cal D}_2(\al)&=&(\mu+2)(\mu+4)p_5(\al)-(15\!+\!15\mu\!+\!3\mu^2)\al p_6(\al)+(24\!+\!27\mu\!+\!6\mu^2)\al^2 p_7(\al)+\cdots\,,\nonumber\\
{\cal D}_3(\al)&=& (\mu+2)(\mu+4)(\mu+6)p_7(\al)+\cdots\, , \label{e34b}
\end{eqnarray}
with $\mu:=b-1$.
Taking account of (\ref{e32a}), we finally obtain the uniform expansion
\begin{theorem}$\!\!\!.$ Let $\epsilon>1$ and the parameters $a$, $b$ and $c$ be real with $b>0$. Further let $0<x<1$. Then we have the uniform expansion valid in the neighbourhood of $\epsilon x=1$ $(\al=0)$
\bee\label{e35}
F(a+\epsilon\la,b;c+\la;x)\sim G(\la)e^{-\la\psi(x,t_s)}\bl\{\frac{W_b(\al\la^\fr)}{\la^{\fr b}} \sum_{k=0}^\infty \frac{{\cal C}_k(\al)}{\la^{k}}+\frac{W_b'(\al\la^\fr)}{\la^{\fr b+\fr}} \sum_{k=0}^\infty \frac{{\cal D}_k(\al)}{\la^{k}}\br\}
\ee
as $\la\to+\infty$, where the first few coefficients ${\cal C}_k(\al)$ and ${\cal D}_k(\al)$ are defined in (\ref{e34a}) and (\ref{e34b}).
\end{theorem}
\vspace{0.2cm}

We remark that the form of the expansion (\ref{e35}) has a similar appearance to that in (\ref{e26}). Indeed this
similarity is even closer since it is shown in Appendix B that the leading coefficients (when formally extended to infinite series) have the values ${\cal C}_0(\al)\to A_0(\al)$, ${\cal C}_1(\al)\to A_1(\al)$, ${\cal D}_0(\al)\to B_0(\al)$ and ${\cal D}_1(\al)\to B_1(\al)$.
The main differences between the expansions in Theorems 1 and 2 are: (i) the coefficients ${\cal C}_k(\al)$, ${\cal D}_k(\al)$ are  easier to compute than $A_k(\al)$, $B_k(\al)$ and (ii) they do not possess a removable singularity when $\al=0$. 
\vspace{0.6cm}

\begin{center}
{\bf 4. \ The expansions at coalescence $\epsilon x=1$}
\end{center}
\setcounter{section}{4}
\setcounter{equation}{0}
\renewcommand{\theequation}{\arabic{section}.\arabic{equation}}
When $\epsilon x=1$ the saddle point is situated at $t=0$ and the quantity $\al=0$. Then for the uniform expansion given in (\ref{e26}) we have from (\ref{e28})--(\ref{e210}) and (\ref{a5}) the coefficients at coalescence given by
\bee\label{e41}
A_0(0)=g_0(0), \quad B_0(0)=g_0'(0),\quad A_1(0)=\frac{1}{2} bg_0''(0), \quad B_1(0)=\frac{1}{6}(1+b) g_0'''(0).
\ee
With $\epsilon x=1+\delta$, $\delta\to0+$, we have 
\[\psi(x,t_s)=\frac{-\epsilon \delta^2}{2(\epsilon-1)}\bl\{1-\frac{2\delta(\epsilon-2)}{3(\epsilon-1)}+O(\delta^2)\br\},\quad\al=\sqrt{\frac{\epsilon}{\epsilon-1}}\,\delta\bl\{1-\frac{\delta(\epsilon-2)}{3(\epsilon-1)}+O(\delta^2)\br\},\]
together with $\psi'(x,0)=-\delta$, $\psi''(x,0)=(\epsilon-1)/\epsilon-2\delta/\epsilon+O(\delta^2)$ and $t_s=\epsilon\delta/(\epsilon-1)+O(\delta^2)$. Then, after some lengthy algebra, we obtain from (\ref{a6}) the values
\bee\label{e41a}
g_0(0)=A_0(0)=\bl(\frac{\epsilon}{\epsilon-1}\br)^{b/2},\qquad g_0'(0)=B_0(0)=\frac{g_0(0)}{\sqrt{\epsilon(\epsilon-1)}}\bl\{a+\frac{1}{3}(b+1)(2\epsilon-1)-\epsilon c\br\}.
\ee
We do not present the details for $g_0''(0)=2! p_2(0)$ and $g_0'''(0)=3! p_3(0)$, where the coefficients $p_r(\al)$ are defined in (\ref{e31}), on account of their algebraic complexity. Their numerical evaluation in specific cases will be discussed in Section 6.

From the properties of the parabolic cylinder function $W_b(\chi)$ in (\ref{e32a}) we have the values
\[W_b(0)=2^{\fr b-1} \frac{\g(\fs b)}{\g(b)},\qquad W_b'(0)=2^{\fr b-\fr} \frac{\g(\fs b+\fs)}{\g(b)}.\]
Then, since $\psi(x,0)=0$, the expansion (\ref{e26}) when $\epsilon x=1$ becomes
\bee\label{e42}
F(a+\epsilon\la,b;c+\la;x)\sim\frac{G(\la)}{2\g(b)}\bl\{\frac{\g(\fs b)}{(\la/2)^{\fr b}}\sum_{k=0}^\infty\frac{A_k(0)}{\la^k}+\frac{\g(\fs b+\fs)}{(\la/2)^{\fr b+\fr}}\sum_{k=0}^\infty\frac{B_k(0)}{\la^k}\br\}
\ee
as $\la\to+\infty$.

For the alternative expansion (\ref{e35}), it is readily seen from (\ref{e33a}) when $\chi=\al\lambda^{1/2}=0$ we have $S_k(0)=(b+k-2) S_{k-2}(0)$, so that
the coefficients $c_k(0)$, $d_k(0)$ take the form
\[c_{2k}(0)=2^k (\fs b)_k,\quad c_{2k+1}(0)=0;\qquad d_{2k}(0)=0,\quad d_{2k+1}(0)=2^k (\fs b+\fs)_k \qquad (k\geq 0).\]
This can also be verified to be the case from Table 1 for the first few $k$-values. Then, when $\epsilon x=1$,
\begin{eqnarray}
F(a+\epsilon\la,b;c+\la;x)&\sim& \frac{G(\la)}{2\g(b)}\bl\{\frac{\g(\fs b)}{(\la/2)^{\fr b}} \sum_{k=0}^\infty\frac{p_{2k}(0)c_{2k}}{\la^k}+\frac{\g(\fs b+\fs)}{(\la/2)^{\fr b+\fr}} \sum_{k=0}^\infty \frac{p_{2k+1}(0) d_{2k+1}}{\la^k}\br\}\nonumber\\
&=&\frac{G(\la)}{2\g(b)}\bl\{\sum_{k=0}^\infty \frac{p_{2k}(0) \g(\fs b\!+\!k)}{(\la/2)^{\fr b+k}}+\sum_{k=0}^\infty\frac{p_{2k+1}(0) \g(\fs b\!+\!k\!+\!\fs)}{(\la/2)^{\fr b+\fr+k}}\br\}\nonumber\\
&=&\frac{G(\la)}{2\g(b)} \sum_{k=0}^\infty \frac{p_k(0) \g(\fs b+\fs k)}{(\la/2)^{\fr b+\fr k}}\qquad (\la\to+\infty),\label{e43}
\end{eqnarray}
where $p_0(0)=g_0(0)$, $p_1(0)=g_0'(0)$ with $g_0(0)$, $g_0'(0)$ given in (\ref{e41a}).

The form of expansion in (\ref{e43}) is, of course, the same as that obtained by application of the method of steepest descents to the integral (\ref{e32}). This procedure produces
\begin{eqnarray*}
F(a+\epsilon\la,b;c+\la;x)&\sim&\frac{G(\la)}{\g(b)} \sum_{k=0}^\infty p_k(0) \int_0^\infty u^{b+k-1} e^{-\fr\la u^2} du\\
&=&\frac{G(\la)}{2\g(b)} \sum_{k=0}^\infty \frac{p_k(0) \g(\fs b+\fs k)}{(\la/2)^{\fr b+\fr k}}\qquad(\la\to+\infty; \ \epsilon x=1),
\end{eqnarray*}
which is equivalent to that obtained in \cite[Section 2]{H}.
\vspace{0.6cm}

\begin{center}
{\bf 5. \ The case $b<0$}
\end{center}
\setcounter{section}{5}
\setcounter{equation}{0}
\renewcommand{\theequation}{\arabic{section}.\arabic{equation}}
The case of negative (non-integer) values of the parameter $b$ is not covered by the expansions in Theorems 1 and 2. When $b<0$ we can employ the contiguous relation \cite[(15.5.11)]{DLMF} expressed in the form
\[{}_2F_1(a+\epsilon\la,b;c+\la;x)=\bl\{1-\epsilon x+\frac{{\bf A}(x)}{\la\!+\!c\!-\!b\!-\!1}\br\}\,{}_2F_1(a+\epsilon\la,b+1;c+\la;x)\]
\bee\label{e51}
\hspace{4cm}+\frac{{\bf B}(x)}{\la\!+\!c\!-\!b\!-\!1}\,{}_2F_1(a+\epsilon\la,b+2;c+\la;x),
\ee
where
\[{\bf A}(x):=(c-b-1)\epsilon x-(a-b-1)x-b-1.\qquad {\bf B}(x)=(1+b)(1-x).\]
Then, if $b\in(-1,0)$, the above relation expresses ${}_2F_1(a+\epsilon\la,b;c+\la;x)$ in terms of similar functions with positive second parameter, from which the uniform expansions in either Theorem 1 or 2 can be separately applied. In the same manner, (\ref{e51}) can be applied twice if $b\in(-2,-1)$ and so on; see also \cite[\S 3.2]{H}.

The non-uniform expansion when $b=-m$, where $m$ is a positive integer, can be obtained from (\ref{e12}). An alternative procedure is to write 
\bee\label{e52}
{}_2F_1(a+\epsilon\la,-m;c+\la;x)=\sum_{r=0}^m(-)^r\bl(\!\!\!\begin{array}{c}m\\r\end{array}\!\!\!\br)(\epsilon x)^r T_r(\la),\qquad T_r(\la):=\prod_{s=0}^{r-1}\frac{(1+\frac{a+s}{\epsilon\la})}{(1+\frac{c+s}{\la})}.
\ee
Observing that
\[\log\,T_r(\la)=\frac{r}{\epsilon\la}\{a-\epsilon c+\fs(r-1)(1-\epsilon)\}\hspace{4cm}\]
\[\hspace{3cm}-\frac{r}{2\epsilon^2\la^2}\{a^2-\epsilon^2 c^2+(a-\epsilon^2c)(r-1)+\f{1}{6}(r-1)(2r-1)(1-\epsilon^2)\}+O(\la^{-3})\]
we find upon exponentiation
\[T_r(\la)=1+\frac{E_0}{\epsilon\la}+\frac{E_1}{2\epsilon^2\la^2}+O(\la^{-3}),\]
where
\[E_0=r(a-\epsilon c)+\fs r(r-1)(1-\epsilon),\qquad E_1=E_0^2-r\{a^2-\epsilon^2 c^2+(a-\epsilon^2c)(r-1)+\f{1}{6}(r-1)(2r-1)(1-\epsilon^2)\}.\]
After straightforward evaluation of the various binomial-type sums that arise, we obtain (provided $\epsilon x\neq 1$)
\bee\label{e53}
{}_2F_1(a+\epsilon\la,-m;c+\la;x)=(1-\epsilon x)^m\bl\{1-\frac{m}{\epsilon\la}\,\sum_{s=1}^2\Omega_{0s}X^s+\frac{m}{2\epsilon^2\la^2}\,\sum_{s=1}^4\Omega_{1s}X^s+O(\la^{-3})\br\}
\ee
as $\la\to+\infty$, where $X:=\epsilon x/(1-\epsilon x)$ and
\[\Omega_{01}=a-\epsilon c,\qquad \Omega_{02}=\fs(m-1)(\epsilon-1),\]
\[\Omega_{11}=a^2-\epsilon^2c^2,\qquad \Omega_{12}=(a-\epsilon c)^2(m-(\epsilon x)^{-1})-(m-1)(a-\epsilon^2c),\]
\[\Omega_{13}=(m-1)(\epsilon-1)\{(a-\epsilon c)(m-2(\epsilon x)^{-1})-\f{1}{6}(1+\epsilon)(2m-1-3(\epsilon x)^{-1})\},\]
\[\Omega_{14}=\f{1}{4}(m-1)(1-\epsilon)^2\{(m-1)(m-4(\epsilon x)^{-1})+2(\epsilon x)^{-2}\}.\]

When $\epsilon x=1$, we find from (\ref{e52}) upon series expansion of $T_r(\la)$ in inverse powers of $\la$ with the aid of {\it Mathematica} the following asymptotic expressions for the function $F_m:={}_2F_1(a+\epsilon\la,-m;c+\la;x)$:
\begin{eqnarray*}
F_1&=&\frac{(\epsilon c-a)}{\epsilon\la}\bl\{1-\frac{c}{\la}+\frac{c^2}{\la^2}+O(\la^{-3})\br\},\\
F_2&=&-\frac{(\epsilon-1)}{\epsilon\la}+\frac{\Upsilon_2}{\epsilon^2\la^2}+O(\la^{-3}),\\
F_3&=&\frac{(\epsilon-1)}{\epsilon^2\la^2}(2+3a+\epsilon(4+3c))-\frac{\Upsilon_3}{\epsilon^3\la^3}+O(\la^{-4}),\\
F_4&=&\frac{3(\epsilon-1)^2}{\epsilon^2\la^2}-\frac{\Upsilon_4}{\epsilon^3\la^3}+O(\la^{-4}),
\end{eqnarray*}
where
\begin{eqnarray*}
\Upsilon_2&=&a(a+1)-2\epsilon a(c+1)+\epsilon(\epsilon-1+\epsilon(3\epsilon-2+c^2),\\
\Upsilon_3&=&a(a+1)(a+2)-3\epsilon a\{a(c+2)+5+4c\}+3\epsilon^2a(c+1)(c+5)\\
&&\hspace{5cm}-\epsilon(c+1)\{6-9\epsilon(c+2)+\epsilon^2(12+11c+c^2)\},\\
\Upsilon_4&=&3a^2+7a+3-2\epsilon(a(7+3c)+9+5c)+\epsilon^2(18+17c+3c^2).
\end{eqnarray*}
The complexity of the coefficients of higher powers of $\la^{-1}$ and higher $m$-values prevents their presentation.
\vspace{0.6cm}

\begin{center}
{\bf 6. \ Numerical results}
\end{center}
\setcounter{section}{6}
\setcounter{equation}{0}
\renewcommand{\theequation}{\arabic{section}.\arabic{equation}}
The coefficients appearing in the expansions in Theorems 1 and 2 depend on the values of the coefficients in the expansion of $g_0(u)$ in ascending powers of $u-\al$. Apart from the leading coefficients, the algebraic complexity of the higher coefficients in general prevents their presentation. However, in specific cases, where the parameters have numerical values, it is perfectly feasible with the help of {\it Mathemtica} to generate many coefficients.

To illustrate, we consider the case $a=c=1$, $b=\f{3}{2}$, $\epsilon=2$ and $x=0.60$, so that from (\ref{e24}) and (\ref{e24b}) the saddle point is at $t_s=\f{1}{3}$ and the parameter $\al=(2\log (25/24))^{1/2}\doteq 0.285734$. Then, from (\ref{e24a}), we have that
\[w\equiv u-\al=\frac{3}{2\surd 2} (t-\frac{1}{3})+\frac{9}{8\surd 2}(t-\frac{1}{3})^2+\frac{135}{128\surd 2}(t-\frac{1}{3})^3+\frac{567}{512\surd 2}(t-\frac{1}{3})^4+\cdots\,,\]
which upon inversion yields
\[t\equiv t(w)=\frac{1}{3}+\frac{2\surd 2}{3}w-\frac{2}{3}w^2+\frac{1}{2\surd 2} w^3-\frac{1}{6} w^4+\frac{29}{288\surd 2} w^5-\frac{1}{36} w^6+\cdots\, .\]
The derivative $dt/dw$ then follows by differentiation of the above expansion. From (\ref{e25a}) and (\ref{e31}),
we finally obtain, upon use of the Series command in {\it Mathematica}, the expansion
\bee\label{e51}
g_0(u)=f(x,t(w)) \bl(\frac{t(w)}{w+\al}\br)^{b-1}\,\frac{dt}{dw}=\sum_{k=0}^\infty p_k(\al) w^k.
\ee
The values of the coefficients $p_k(\al)$ for $k\leq 7$ are displayed in Table 2, from which those of $g_0^{(k)}(\al)$ can be deduced since $p_k(\al)=g_0^{(k)}(\al)/k!$. 
\begin{table}[th]
\caption{\footnotesize{Values of the coefficients $p_k(\al)$ for $k\leq 7$ when $a=c=1$, $b=3/2$, $\epsilon=2$ and (i) $x=0.60$\ (ii) $x=0.50$}} \label{t2}
\begin{center}
\begin{tabular}{|l|l|l||l|l|l|}
\hline
%&&\\[-0.3cm]
\mcol{1}{|c|}{} & \mcol{1}{c|}{$\epsilon x=1.20$} & \mcol{1}{c||}{$\epsilon x=1.00$} &  \mcol{1}{c|}{} & \mcol{1}{c|}{$\epsilon x=1.20$} & \mcol{1}{c|}{$\epsilon x=1.00$}\\
\mcol{1}{|c|}{$k$} & \mcol{1}{c|}{$p_k(\al)$} & \mcol{1}{c||}{$p_k(0)$} & \mcol{1}{c|}{$k$} & \mcol{1}{c|}{$p_k(\al)$} & \mcol{1}{c|}{$p_k(0)$}\\
\hline
&&&&&\\[-0.3cm]
0 & 2.3384546881 & 1.6817928305 & 4 & 0.1073018863 & 0.1572235514\\
1 & 2.5221370759 & 1.7838106725 & 5 & 0.0348400522 & 0.0513462836\\
2 & 1.2129759615 & 0.8408964153 & 6 & 0.0098542549 & 0.0145775887\\
3 & 0.4345143534 & 0.2973017788 & 7 & 0.0027097818 & 0.0040044129\\
[.1cm]\hline\end{tabular}
\end{center}
\end{table} 

The expansion given in Theorem 1 requires the calculation of the coefficients $A_r(\al)$, $B_r(\al)$ ($r=1, 2$) defined in (\ref{e28})--(\ref{e210}). With the parameter values mentioned above, we find the values:
\[A_0(\al)=2.3384546881,\quad B_0(\al)=2.2076726817,\]
\[A_1(\al)=1.7632504455,\quad B_1(\al)=0.9274745592.\] 
These values of the coefficients have been used in constructing the entry corresponding to $x=0.60$ in Table 3, which shows the absolute relative error\footnote{In the tables we write the values as $x(y)$ instead of $x\times 10^y$.} in the computation of ${}_2F_1(a+\epsilon\la,b;c+\la;x)$ using the expansion in (\ref{e26}) with $k\leq 1$ for a range of $x$ and $\la$ values. The value $x=0.50$ in Table 3 corresponds to $\epsilon x=1$ ($\al=0$); the values of the coefficients in this case are obtained in a similar manner and, from (\ref{e41}), are found to be:
\[A_0(0)=1.6817928305,\quad B_0(0)=1.7838106725,\]
\[A_1(0)=1.2613446229,\quad B_1(0)=0.7432544469.\]

\begin{table}[th]
\caption{\footnotesize{Values of the absolute relative error in the computation of ${}_2F_1(a+\epsilon\la,b;c+\la;x)$ using the expansion (\ref{e26}) with $k\leq 1$ when $a=c=1$, $b=3/2$, $\epsilon=2$ for different $x$ and $\la$.}} \label{t3}
\begin{center}
\begin{tabular}{|l|lllll|}
\hline
&&&&&\\[-0.3cm]
\mcol{1}{|c|}{} & \mcol{5}{c|}{$\la$}\\
\mcol{1}{|c|}{$x$} & \mcol{1}{c}{10} & \mcol{1}{c}{20} & \mcol{1}{c}{50} & \mcol{1}{c}{100} & \mcol{1}{c|}{200}\\
\hline
&&&&&\\[-0.3cm]
0.70 & $3.215(-03)$ & $8.227(-04)$ & $1.335(-04)$ & $3.355(-05)$ & $8.407(-06)$\\
0.60 & $3.128(-03)$ & $8.058(-04)$ & $1.321(-04)$ & $3.331(-05)$ & $8.363(-06)$\\
0.55 & $3.070(-03)$ & $7.915(-04)$ & $1.302(-04)$ & $3.300(-05)$ & $8.317(-06)$\\
0.50 & $3.004(-03)$ & $7.733(-04)$ & $1.272(-04)$ & $3.228(-05)$ & $8.155(-06)$\\
0.45 & $2.933(-03)$ & $7.526(-04)$ & $1.234(-04)$ & $3.122(-05)$ & $7.869(-06)$\\
0.40 & $2.860(-03)$ & $7.309(-04)$ & $1.192(-04)$ & $3.007(-05)$ & $7.559(-06)$\\
0.30 & $2.711(-03)$ & $6.873(-04)$ & $1.112(-04)$ & $2.792(-05)$ & $6.997(-06)$\\
[.1cm]\hline\end{tabular}
\end{center}
\end{table} 

For the expansion in (\ref{e35}), the coefficients ${\cal C}_k(\al)$ and ${\cal D}_k(\al)$ in (\ref{e34a}) and (\ref{e34b}) are computed using the values of $p_k(\al)$ obtained in (\ref{e51}). The absolute relative error in the computation of ${}_2F_1(a+\epsilon\la,b;c+\la;x)$ using the expansion (\ref{e35}) with $k\leq 3$ is shown in Table 4 for a range of $x$ and $\la$ values.

\begin{table}[th]
\caption{\footnotesize{Values of the absolute relative error in the computation of ${}_2F_1(a+\epsilon\la,b;c+\la;x)$ using the expansion (\ref{e35}) with $k\leq 3$ when $a=c=1$, $b=3/2$, $\epsilon=2$ for different $x$ and $\la$.}} \label{t3}
\begin{center}
\begin{tabular}{|l|lllll|}
\hline
&&&&&\\[-0.3cm]
\mcol{1}{|c|}{} & \mcol{5}{c|}{$\la$}\\
\mcol{1}{|c|}{$x$} & \mcol{1}{c}{10} & \mcol{1}{c}{20} & \mcol{1}{c}{50} & \mcol{1}{c}{100} & \mcol{1}{c|}{200}\\
\hline
&&&&&\\[-0.3cm]
0.70 & $1.956(-05)$ & $1.401(-06)$ & $3.499(-08)$ & $2.124(-09)$ & $1.317(-10)$\\
0.60 & $2.350(-05)$ & $1.834(-06)$ & $6.043(-08)$ & $4.132(-09)$ & $2.435(-10)$\\
0.55 & $2.755(-05)$ & $2.248(-06)$ & $7.797(-08)$ & $5.978(-09)$ & $4.478(-10)$\\
0.50 & $3.471(-05)$ & $3.166(-06)$ & $1.320(-07)$ & $1.185(-08)$ & $1.059(-09)$\\
0.45 & $4.789(-05)$ & $5.387(-06)$ & $3.615(-07)$ & $6.042(-08)$ & $1.443(-08)$\\
0.40 & $7.366(-05)$ & $1.143(-05)$ & $1.616(-06)$ & $6.117(-07)$ & $3.461(-07)$\\
0.30 & $2.492(-04)$ & $8.563(-05)$ & $3.927(-05)$ & $2.972(-05)$ & $2.581(-05)$\\
[.1cm]\hline\end{tabular}
\end{center}
\end{table} 

\newpage
\vspace{0.6cm}

\begin{center}
{\bf Appendix A: Derivation of the coefficients $A_1(\al)$ and $B_1(\al)$}
\end{center}
\setcounter{section}{1}
\setcounter{equation}{0}
\renewcommand{\theequation}{\Alph{section}.\arabic{equation}}
From (\ref{e25b}) the function $g_0(u)$ is specified by
\[g_0(u)=A_0(\al)+B_0(\al)(u-\al)+u(u-\al) G_0(u),\]
where $A_0(\al)=g_0(\al)$, $B_0(\al)=\{g_0(\al)-g_0(0)\}/\al$. Then we have
\[G_0(u)=\frac{g_0(u)-A_0(\al)-B_0(\al)(u-\al)}{u(u-\al)}
=\frac{\al g_0(u)-ug_0(\al)+(u-\al)g_0(0)}{\al u(u-\al)}.\]
Following \cite[\S 22.2]{T}, we may write this as the Cauchy integral
\[G_0(u)=\frac{1}{2\pi i}\int_C \frac{g_0(s)}{s(s-\al)(s-u)}\,du,\]
where $C$ denotes a closed contour in the domain of analyticity of $g_0(s)$ containing the points $0$, $\al$ and $u$ in its interior.

From (\ref{e25c}), the next function in the sequence is
\begin{eqnarray*}
G_1(u)&=&bG_0(u)+uG_0'(u)\\
&=&\frac{1}{2\pi i}\int_C \frac{b(s-u)+u}{s(s-\al)(s-u)^2}\,g_0(s)\,ds.
\end{eqnarray*}
Evaluation of the residues then yields
\[g_1(u)=\frac{(b-1)}{\al u}\,g_0(0)+\frac{\al b+(1-b)u}{\al(u-\al)^2}\,g_0(\al)\]
\bee\label{a1}
+\frac{1}{(u-\al)^2}\bl\{\bl(\frac{\al(1-b)}{u}+b-2\br)\,g_0(u)+(u-\al)g_0'(u)\br\}.
\ee
Upon insertion of the expansion $g_0(u)=g_0(0)+ug_0'(0)+\ldots\,$ and similarly for $g_0'(u)$, we obtain after some routine algebra
\bee\label{a2}
g_1(0)=\frac{b}{\al^2}\{g_0(\al)-g_0(0)\}-\frac{b}{\al}\,g_0'(0)=\frac{b}{\al}\{B_0(\al)-g_0'(0)\}.
\ee

A similar procedure using the expansion $g_0(u)=g_0(\al)+(u-\al)g_0'(\al)+\fs (u-\al)^2 g_0''(\al)+\ldots$ leads to 
\[g_1(\al)=\frac{(1-b)}{\al^2}\{g_0(\al)-g_0(0)\}-\frac{(1-b)}{\al}\,g_0'(\al)+\frac{1}{2} g_0''(\al)\]
\bee\label{a3}
=\frac{(1-b)}{\al}\{B_0(\al)-g_0'(\al)\}+\frac{1}{2}\,g_0''(\al).
\ee
Then, from (\ref{e29}), the coefficients $A_1(\al)$ and $B_1(\al)$ are given by
\bee\label{a4}
A_1(\al)=g_1(\al),\qquad B_1(\al)=\frac{1}{\al}\{g_1(\al)-g_1(0)\}.
\ee
\vspace{0.1cm}

\noindent{\bf Values of the coefficients at coalescence $\al=0$.}\ \ 
From (\ref{e28}), it follows that
\[B_0(\al)=g_0'(0)+\frac{1}{2}\al g_0''(0)+\frac{1}{6}\al^2 g_0'''(0)+O(\al^3),\]
so that from (\ref{a2}) and (\ref{a3}) we find after some algebra the expansions
\[g_1(0)=\frac{1}{2}b g_0''(0)+\frac{1}{6}\al b g_0'''(0)+O(\al^2),
\qquad g_1(\al)=\frac{1}{2}b g_0''(0)+\frac{1}{6}\al (1+2b) g_0'''(0)+O(\al^2).\]
Then, from (\ref{a4}) the coefficients at coalescence ($\al=0$) are given by
\bee\label{a5}
A_1(0)=\frac{1}{2}bg_0''(0),\qquad B_1(0)=\frac{1}{6}(1+b) g_0'''(0).
\ee

To determine the value of $g_0'(0)$ at coalescence, we note that from (\ref{e24a})
\[\psi(x,t)=t\psi'(x,0)+\frac{1}{2}t^2 \psi''(x,0)+\cdots=-\al u+\frac{1}{2}u^2\]
and hence, upon inversion,
\[t=-\frac{\al}{\psi'(x,0)}u+\bl(1-\frac{\al^2\psi''(x,0)}{\psi'(x,0)^2}\br)\,\frac{u^2}{2\psi'(x,0)}+\cdots\ .\]
Substitution of this last expansion in $g_0(u)$ in (\ref{e25a}) then produces
\[g_0(u)=g_0(0)\bl(1+\bl\{\frac{(b+1)}{2\al}\bl(\frac{\al^2\psi''(x,0)}{\psi'(x,0)^2}-1\br)-\frac{\al f'(x,0)}{\psi'(x,0)}\br\}u+\cdots\br)\]
to yield the value 
\bee\label{a6}
g_0'(0)=g_0(0)\bl\{\frac{(b+1)}{2\al}\bl(\frac{\al^2\psi''(x,0)}{\psi'(x,0)^2}-1\br)-\frac{\al f'(x,0)}{\psi'(x,0)}\br\},
\ee
where
\[g_0(0)=\bl(-\frac{\al}{\psi'(x,0)}\br)^b, \qquad f'(x,0)=ax+b+1-c.\]
\vspace{0.6cm}

\begin{center}
{\bf Appendix B: Limiting values of the coefficients ${\cal C}_k(\al)$ and ${\cal D}_k(\al)$ when $k=0, 1$}
\end{center}
\setcounter{section}{2}
\setcounter{equation}{0}
\renewcommand{\theequation}{\Alph{section}.\arabic{equation}}
If we formally extend the leading coefficients ${\cal C}_1(\al)$, ${\cal D}_0(\al)$ and ${\cal D}_1(\al)$ in Theorem 2 to infinite sums (which we denote by the addition of a circumflex) we have
\[{\hat{\cal C}}_1(\al)=p_2(\al)+\frac{\mu}{\al^2} \sum_{k\geq 2}(-)^kp_k(\al) \al^k,\quad {\hat{\cal D}}_0(\al)=-\frac{1}{\al}\sum_{k\geq1}(-)^kp_k(\al)\al^k,\]
\[{\hat{\cal D}}_1(\al)=-\frac{1}{\al^3}\sum_{k\geq 3}(-)^k p_k(\al) \al^k \{k(1+\mu)-(1+2\mu)\},\]
where we recall that $\mu:=b-1$. From (\ref{e31}), we have the expansions
\[g_0(0)=\sum_{k\geq0}(-)^kp_k(\al) \al^k,\qquad g_0'(0)=-\frac{1}{\al}\sum_{k\geq 1}(-)^kp_k(\al) \al^k,\qquad p_k(\al)=\frac{g_0^{(k)}(\al)}{k!},\]
and from (\ref{e28}) the values of the leading coefficients
\[A_0(\al)=g_0(\al),\qquad B_0(\al)=\frac{1}{\al}\{g_0(\al)-g_0(0)\}.\]

First, we have
\bee\label{b1}
{\hat{\cal C}}_0(\al)=p_0(\al)=g_0(\al)=A_0(\al).
\ee 
Then it follows that
\begin{eqnarray}
{\hat{\cal C}}_1(\al)&=&p_2(\al)+\frac{\mu}{\al^2}\bl\{\sum_{k\geq0}(-1)^kp_k(\al) \al^k-p_0(\al)+\al p_1(\al)\br\}\nonumber\\
&=&\frac{\mu}{\al^2}\{g_0(0)-g_0(\al)+\al g_0'(\al)\}+\fs g_0''(\al)\nonumber\\
&=&\frac{\mu}{\al}\{g_0'(\al)-B_0(\al)\}+\fs g_0''(\al)=g_1(\al)\nonumber\\
&=&A_1(\al)\label{b2}
\end{eqnarray}
by (\ref{e29}).

Similarly, we obtain
\bee\label{b3}
{\hat{\cal D}}_0(\al)=-\frac{1}{\al}\bl\{\sum_{k\geq0}(-)^kp_k(\al) \al^k+p_0(\al)\br\}=\frac{1}{\al}\{g_0(\al)-g_0(0)\}=B_0(\al)
\ee
and
\begin{eqnarray}
{\hat{\cal D}}_1(\al)&=&-\frac{1}{\al^3}\sum_{k\geq0}(-)^kp_k(\al) \al^k\{k(1+\mu)-(1+2\mu)\}\nonumber\\
&&\hspace{5cm} +\frac{1}{\al^3}\{-(1+2\mu)p_0(\al)+\mu\al p_1(\al)+p_2(\al) \al^2\}\nonumber\\
&=&\frac{1}{\al^3}\bl\{(1+2\mu)\{g_0(0)-g_0(\al)\}+(1+\mu)\al g_0'(0)+\mu\al g_0'(\al)+\fs\al^2 g_0''(\al)\br\}\nonumber\\
&=&-\frac{(1+2\mu)}{\al^2} B_0(\al)+\frac{(1+\mu)}{\al^2} g_0'(0)+\frac{\mu g_0'(\al)}{\al^2}+\frac{1}{2\al} g_0''(\al)\nonumber\\
&=&-\frac{(1+\mu)}{\al^2}\{B_0(\al)-g_0'(0)\}+\frac{\mu}{\al^2}\{g_0'(\al)-B_0(\al)\}+\frac{1}{2\al} g_0''(\al)\nonumber\\
&=&\frac{1}{\al}\{g_1(\al)-g_1(0)\}=B_1(\al)\label{b4}
\end{eqnarray}
by (\ref{e29}) and (\ref{e210}).

\end{document}